\documentclass[runningheads]{llncs}

\usepackage[T1]{fontenc}
\usepackage{graphicx}
\usepackage{amsmath}
\usepackage{amssymb}
\usepackage{hyperref}
\usepackage{algorithm,algpseudocode} 
\usepackage{subcaption}
\usepackage{xcolor}
\usepackage{lineno}
\usepackage{arydshln}
\usepackage{pifont}   
\usepackage{colortbl} 

\definecolor{purple}{RGB}{160,32,240}
\definecolor{orange}{RGB}{255,165,0}
\definecolor{green}{RGB}{0,128,0}

\newcommand{\cmark}{\textcolor{green}{\checkmark}}
\newcommand{\xmark}{\textcolor{purple}{\ding{55}}}
\newcommand{\approxmark}{\textcolor{orange}{$\sim$}}

\title{Empirical and computer-aided robustness analysis of long-step and accelerated methods\\in smooth convex optimization}
\titlerunning{Robustness analysis of long-step and accelerated optimization methods}
\author{Pierre Vernimmen \and François Glineur}
\institute{Mathematical Engineering Department, UCLouvain, Louvain-la-Neuve, Belgium\newline
\email{\{pierre.vernimmen,francois.glineur\}@uclouvain.be}}

\begin{document}

\maketitle
\begin{abstract} This work assesses both empirically and theoretically, using the performance estimation methodology, how robust different first-order optimization methods are when subject to relative inexactness in their gradient computations. Relative inexactness occurs, for example, when compressing the gradient using fewer bits of information, which happens when dealing with large-scale problems on GPUs. Three major families of methods are analyzed: constant step gradient descent, long-step methods, and accelerated methods. The latter two are first shown to be theoretically not robust to inexactness. Then, a semi-heuristic shortening factor is introduced to improve their theoretical guarantees. All methods are subsequently tested on a concrete inexact problem, with two different types of relative inexactness, and it is observed that both accelerated methods are much more robust than expected, and that the shortening factor significantly helps the long-step methods. In the end, all shortened methods appear to be promising, even in this inexact setting.

\end{abstract}
\vspace{-0.2cm}
\section{Introduction and Related Work}
\vspace{-0.2cm}
\subsection{Inexact gradient methods} \label{sec::intro}

We first recall that a function \(f : \mathbb{R}^d \to \mathbb{R}\) is said to be \emph{convex} if, for all \(x, y \in \mathbb{R}^d\) and any \(\alpha \in [0, 1]\), the following inequality holds:
\[
f(\alpha x + (1 - \alpha)y) \le \alpha f(x) + (1 - \alpha)f(y).
\]
This property ensures that the line segment between any two points on the graph of \(f\) lies above the graph itself, which implies that local minima are also global minima.

Moreover, we say that \(f\) is \(L\)-\emph{smooth} (or equivalently, that its gradient is Lipschitz continuous) if there exists a constant \(L > 0\) such that for all \(x, y \in \mathbb{R}^d\),
\[
\|\nabla f(x) - \nabla f(y)\| \le L \|x - y\|.
\]
This condition guarantees that the gradient of \(f\) does not change too rapidly, which is essential for controlling the behavior and convergence rate of gradient-based methods.

In this work, we consider several first-order methods applied to $L$-smooth convex functions. Our goal is to study how convergence is affected when the gradient used at each iteration is computed inexactly. In particular, we will study the convergence rate of the methods after $N$ iterations in (squared) gradient norm, which is the smallest quantity  $\tau_N$ 
such that the following holds (note the use of a $\min$ to deal with non-monotone methods):
\begin{equation}\label{eq::def_rate}
        \tfrac{1}{L}\min_{k \in \{0,...,N\}} \|\nabla f(x_k)\|^2 \leq \tau_N (f(x_0)-f(x_*)).
\end{equation}

Inexact gradients occur in a large variety of situations, such as the use of floating point computations with limited accuracy (see Section \ref{sec::truncated_mantissa}), dependence on data that is only known approximately, or more generally when facing the time-accuracy tradeoff that is almost always present when the objective function (and its gradient) is obtained through another iterative procedure (e.g.\@ a simulation or another optimization process). Stochasticity can also be viewed as a particular case of inexactness.

Some first-order methods relying on an inexact gradient have been studied before theoretically, using several distinct notions of inexactness. The approximate gradient introduced in \cite{doi:10.1137/060676386} is assumed to differ from the true gradient by some error whose norm is bounded. In \cite{devolder2014first}, another notion of inexact gradient is developed, based on the maximal error incurred by the corresponding quadratic upper bound. In these two cases, the approximation error is not directly related to the scale of the gradient, i.e.\@ error is measured in an \emph{absolute} manner.

In this work, we focus on a \emph{relative} notion of inexactness, where the norm of the difference between the true gradient and its approximate value is bounded by a certain fraction \footnote{Values of the fraction $\delta \ge 1$ do not make sense, as they allow $d_k = 0$, leading to methods making no progress.} of the gradient norm. The inexact gradient at $x_k$ is written $d_k$, and satisfies  the following inequality: 
\begin{equation}\label{rel_inexact}
    \|d_k - \nabla f(x_k)\| \leq \delta \|\nabla f(x_k)\|\text{, with relative inexactness level } \delta \in[0,1)
\end{equation}

This notion captures the idea that the precision of the gradient approximation should scale with the gradient norm, which is especially useful in optimization contexts where the gradient becomes small near an optimum. It was considered in \cite{de2020worst} to analyze the effect of inexactness when dealing with smooth, strongly convex functions.

Our results show that in the smooth convex setting, even with substantial inexactness, classical first-order methods can still converge at a good rate, which is important for real-world optimization problems where exact gradients are often unavailable. We also illustrate how truncating the mantissa of the gradient components naturally leads to this type of relative inexactness (Section \ref{sec::truncated_mantissa}). More specifically, we analyze two types of gradient-type methods, namely fixed step-size gradient descent (Algorithm \ref{algo::inexact_gradient}) and momentum-based accelerated gradient descent (Algorithm \ref{algo::accelerated_inexact_gradient}, from \cite{nesterov1983method}), both using the worst-case theoretical PEP methodology (Sections \ref{sec::intro_PEP} and \ref{sec::theory}) and on a smooth convex logistic classification problem (Section \ref{sec::experiments}).

\vspace{-0.5cm}
\begin{algorithm}[H]
\caption{Inexact Gradient Descent with Relative Inexactness $\delta$}
\label{algo::inexact_gradient}
\begin{algorithmic}[1]
\State Given an $L$-smooth convex function $f$, an initial iterate $x_0$, a schedule of step sizes $h_k$, a number of iterations $N$, and an inexactness parameter $\delta \in[0,1)$ 
\State $k \gets 0$
\While{$k < N$}
    \State Compute an approximate gradient $d_k$ at $x_k$ with relative inexactness:
    \begin{equation*}
    \|d_k - \nabla f(x_k)\| \leq \delta \|\nabla f(x_k)\|
    \end{equation*}
    \vspace{-0.3cm}

    \State Update parameter:
    \vspace{-0.35cm}
    \begin{equation*}\label{eq::update_GD}
    x_{k+1} = x_k - \frac{h_k}{L} d_k
    \end{equation*}

    \State $k \gets k + 1$
\EndWhile
\State \Return Iterate $x_k$ with the smallest gradient norm 
\end{algorithmic}
\end{algorithm}
\vspace{-1.3cm}

\begin{algorithm}
\caption{Inexact Fast Gradient Descent with Relative Inexactness $\delta$}
\label{algo::accelerated_inexact_gradient}
\begin{algorithmic}[1]
\State Given an $L$-smooth convex function $f$, an initial iterate $x_0$, a step size $h$, a number of iterations $N$, and an inexactness parameter $\delta \in[0,1)$
\State $k \gets 0$
\State $y_0 = x_0$
\While{$k < N$}
    \State Compute an approximate gradient $d_k$ at $y_k$ with relative inexactness:
    \begin{equation*}
    \|d_k - \nabla f(y_k)\| \leq \delta \|\nabla f(y_k)\|
    \end{equation*}
    \vspace{-0.5cm}

    \State Update parameters:
    \vspace{-0.5cm}
    \begin{equation*}
    \begin{aligned}
        x_{k+1} & = y_k - \frac{h}{L} d_k\\
        y_{k+1} &= x_{k+1} + \frac{k-1}{k+2} (x_{k+1} - x_{k})
    \end{aligned}
    \end{equation*}
    \State $k \gets k + 1$
\EndWhile
\State \Return Iterate $x_k$ with the smallest gradient norm 
\end{algorithmic}
\end{algorithm}
\vspace{-1.2cm}

\subsection{An application of inexact gradient: compressed gradient descent}\label{sec::truncated_mantissa}

One way relative inexactness naturally arises is through compression techniques, which deliberately reduce gradient precision to save memory, such as in Compressed Gradient Descent, where each of the components of the gradient is compressed, as can be seen hereunder.
Computers store numbers with limited precision. In 32-bit floating point, a  number $a$ is represented with a sign, a mantissa $b$, and an exponent $e$: \[
        a={\underbrace{\pm}_{\text{1 bit}} \times \underbrace{b}_{\text{23 bits}} \times \underbrace{2^{e - 127}}_{\text{8 bits}}}, \quad b \in [1,2)
\]

If we approximate $a$ by $\tilde{a}$, using, for example, $0$ bit of information instead of $23$ for the mantissa (this amounts to fixing $b=1$), the maximal possible relative error in this context of compressed gradient descent is as follows:

\begin{equation*}
    \delta_{\text{compressed}} = \frac{|a-\tilde{a}|}{|a|}=\frac{|\pm b2^{e-127}-\pm 2^{e-127}|}{|\pm b2^{e-127}|}=\frac{b-1}{b}\in \left[0;\frac{1}{2}\right[< \frac{1}{2}
\end{equation*}
This means we have a relative error of $\frac{1}{2}$ at most for each $a$. Furthermore, it can be shown that if each gradient component is $\delta$-relatively inexact, then the whole gradient also is. In the example, when each gradient component’s mantissa is fixed to 1 (0 bit) reduces storage from $32d$ to $9d$ bits per gradient vector, where $d$ is the number of components. If Algorithms \ref{algo::inexact_gradient} and \ref{algo::accelerated_inexact_gradient} maintain strong performance compared to a version with more accurate gradients, the potential memory savings, reducing storage by more than a factor of three, can be substantial, particularly for high-dimensional optimization problems. More generally, it can be shown that keeping $n_\text{bit}$ bits in each mantissa leads to an inexact gradient with $\delta_{\text{compressed}}=\left(\frac{1}{2}\right)^{n_\text{bit}+1}$. Such limited-precision floating-point computations are now frequently offered in modern hardware, such as in graphics processing unit (GPUs), e.g. with the FP64, FP32 and FP16 modes. \vspace*{-.5cm}

\begin{figure}[H]
    \centering
    \includegraphics[width=0.4\linewidth]{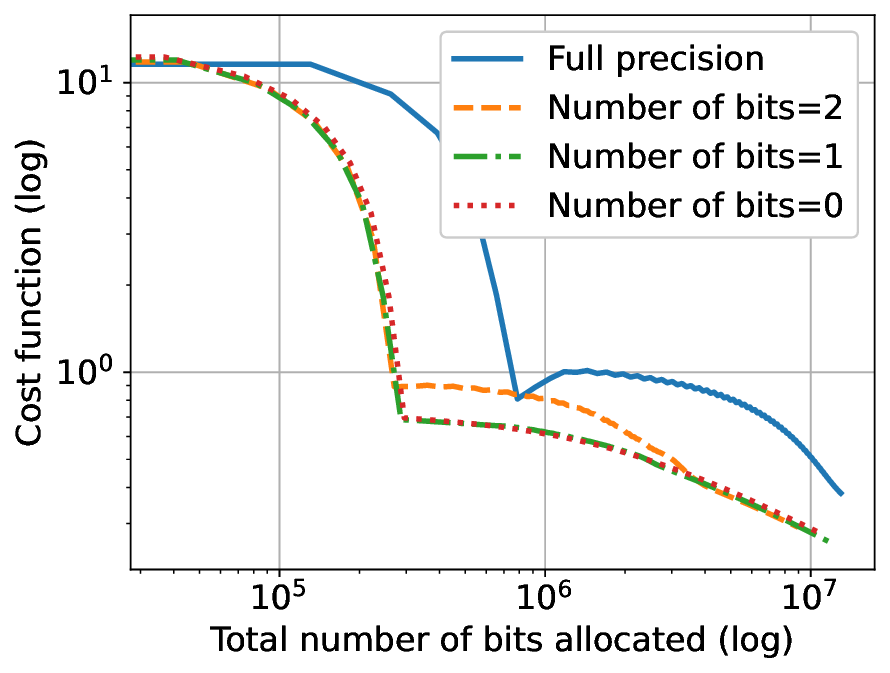}
    \caption{\centering Convergence of the cost function in Gradient Descent (full precision) and Compressed Gradient Descent (with $2$-, $1$- and $0$-bit mantissa), as a function of the total number of bits allocated throughout the execution of the algorithm}
    \label{Experiment_mantissa}
\end{figure} \vspace*{-.5cm}

A preliminary numerical experiment was conducted on a logistic regression problem of classifying hand pictures showing numbers $0$ and $1$ \cite{arda_mavi_2017}. In Figure \ref{Experiment_mantissa}, we plot the cost function as a function of the total number of bits allocated in memory throughout the execution of both Gradient Descent (full precision) and Compressed Gradient Descent\footnote{Total number of bits = [number of iterations performed] $\times$ [ 9 + bits used for the mantissa] $\times$ [$\#$ elements in the gradient]}. It can be seen that using Compressed Gradient Descent improves the method's performance compared to classical Gradient Descent, thus showing the usefulness of analyzing such relatively inexact methods. In experiments (Section \ref{sec::exp_compressed}), we show that such a compression leading to significant memory savings still maintains a competitive performance in practical problems.

\subsection{Description of the considered first-order methods }\label{sec::schedules}

In this work, we analyze several first-order methods with distinct step size schedules, all of which were originally designed for the exact case ($\delta = 0$), and are first evaluated in this context to establish a baseline. Following this, they are modified and adjusted to handle the specific requirements of relative inexactness. Evaluating them in the inexact setting allows us to measure their robustness and effectiveness under relative inexactness. Each step size is normalized by the Lipschitz constant $L$ of the gradient \cite{nesterov2013introductory}.

\begin{itemize}
\item \textbf{Constant Step Size:} We use the simple choice where $h_n = 1.5 \quad \forall n$. This is optimal for one step of Algorithm \ref{algo::inexact_gradient} without inexactness ($\delta=0$) \cite{taylor2017smooth}. The (squared) gradient norm of the iterates can be shown to converge in $\mathcal{O}(\frac{1}{N})$ after $N$ iterations in the exact case (using the $\tau$ convergence measure described above). \medskip

\item \textbf{Dynamic Step Size:} From Teboulle and Vaisbourd \cite{teboulle2023elementary}, this fixed, constant stepsize schedule is defined as:
\begin{equation*}
h_0 = \sqrt{2}, \quad h_n = \frac{-H_{n-1} + \sqrt{H_{n-1}^2 + 8(H_{n-1} + 1)}}{2},
\end{equation*}
where $H_{n-1} := \sum_{i=0}^{n-1} h_i$, and also converges in $\mathcal{O}(\frac{1}{N})$. Its worst-case performance is close to the one achieved by an optimal choice of a constant stepsize, but without requiring knowing the number of performed iterations in advance \cite{teboulle2023elementary}. As $N \to \infty$, step sizes tend towards $2$, which is known to be the optimal constant step size in the exact case when performing $N \to \infty$ iterations \cite{taylor2017smooth}. \medskip

\item \textbf{Silver Step Size Schedule:} this long-step schedule was recently proposed by Altschuler and Parrilo \cite{altschuler2024acceleration}. For $n=2^k-1$, it can be defined recursively:
\begin{equation*}
\begin{array}{rl}
h_1 &:= [\sqrt{2}] \\
h_{1:2n+1} &:= [h_{1:n}, 1 + \rho^{k-1}, h_{1:n}],
\end{array}
\end{equation*}
where $\rho := 1 + \sqrt{2}$ is the \textit{silver ratio}. This schedule relies on performing very long steps from time to time, interspersed with a lot of smaller step sizes. It accelerates convergence to $\mathcal{O}(\frac{1}{N^{1.27}})$ in the exact case, but only guarantees this worst-case performance after given numbers of iterations (those satisfying $n=2^k-1$). \medskip

\item \textbf{Fast Gradient Method (FGM):} This is the classical fast gradient method of Nesterov (Algorithm \ref{algo::accelerated_inexact_gradient}) with $h=1$, which achieves the optimal $\mathcal{O}(\frac{1}{N^2})$ convergence rate in objective function values in the exact case ($\delta=0$). It however requires keeping an additional previous iterate in memory and will serve as a baseline to determine how relative inexactness impacts accelerated methods.
\end{itemize} 

\subsection{Introduction to the PEP methodology}\label{sec::intro_PEP}
All algorithms described in Section \ref{sec::schedules} are first analyzed using the so-called Performance Estimation methodology \cite{drori2014performance}, which allows the automated computation of exact worst-case performance guarantees for a very large class of first-order optimization algorithms. The approach relies on the observation that, by definition, the worst-case behavior of a first-order black-box optimization algorithm is by itself an optimization problem that consists of finding the worst (maximal) convergence rate (recall its definition in \eqref{eq::def_rate}) over all possible inputs to the algorithm. In conclusion, $x_N$ being the output of an algorithm $\mathcal{A}$ after making $N$ calls to the oracle $\mathcal{O}_f$, we look after the worst possible function $f$ (with iterates $x_0,\ldots,x_N$) within a given function class $\mathcal{F}$, equivalent to solving the following Performance Estimation Problem (PEP) \cite{taylor2017smooth}:\vspace*{-.1cm}

\begin{equation}
    \begin{aligned}\label{eq::PEP_opti_problem}
    \tau_N = \max_{f\in \mathcal{F}, \{ x_0,\ldots,x_N,x_\ast \} \subset \mathbb{R}^d} \quad & \min_{k \in \{0,...,N\}} \|\nabla f(x_k)\|^2\\
    \textrm{s.t.} 
    \quad &x_{1}, \cdots, x_N \textrm{ are generated from } x_0 \textrm{ by algorithm }\mathcal{A} \\
    & x_\ast \in \arg\min_{x \in \mathbb{R}^d} f(x), f(x_0) - f(x_\ast) \leq 1.
    \end{aligned}
    \tag{PEP}
\end{equation}

Solving the optimization problem (\ref{eq::PEP_opti_problem}) gives the tight worst-case performance of a first-order algorithm $\mathcal{A}$ on the class of functions $\mathcal{F}$, such as the class of smooth convex functions. This worst-case performance is said to be tight, as it is exactly equal to the worst-case rate. The last constraint bounds the difference in objective function value between the initial point $x_0$ and the optimal solution $x_\ast$, as it is well known that in most situations, the performance of a first-order method cannot be sensibly assessed without such a constraint. The distance between $x_0$ and $x_\ast$ can also be considered, leading to potentially different rates. Still, in this work, we focus on the objective function difference, leading to the convergence rate $\tau_N$ introduced earlier.

To solve problem (\ref{eq::PEP_opti_problem}), which contains an infinite-dimensional variable (the function $f$), it must first be reformulated using interpolation conditions characterizing the class of functions $\mathcal{F}$, as described in \cite{taylor2017smooth}. Inexactness of the gradient can also be dealt with in this setting by writing condition \eqref{rel_inexact} as an additional constraint \cite{de2020worst,taylor2017exact}. This leads to a tractable semidefinite program to compute the tight convergence rate $\tau_N$ and obtain insights about worst-case behaviors.


\vspace{-0.2cm}
\section{Worst-case Results}\label{sec::theory}
\vspace{-0.2cm}

The above-described PEP methodology allows us to get the numerical worst-case performance of a given method with fixed parameters. In practice we used the PEPit toolbox \cite{goujaud2024pepit} with the MOSEK solver \cite{mosek}. Given that the schedule of step sizes presented in Section \ref{sec::schedules} is already fixed,  fixing different values of inexactness $\delta$ allows us to get the theoretical convergence rate $\tau_N$.

\begin{figure}[H]
    \centering
    \begin{subfigure}{0.49\textwidth}
    \centering
        \includegraphics[scale=0.35]{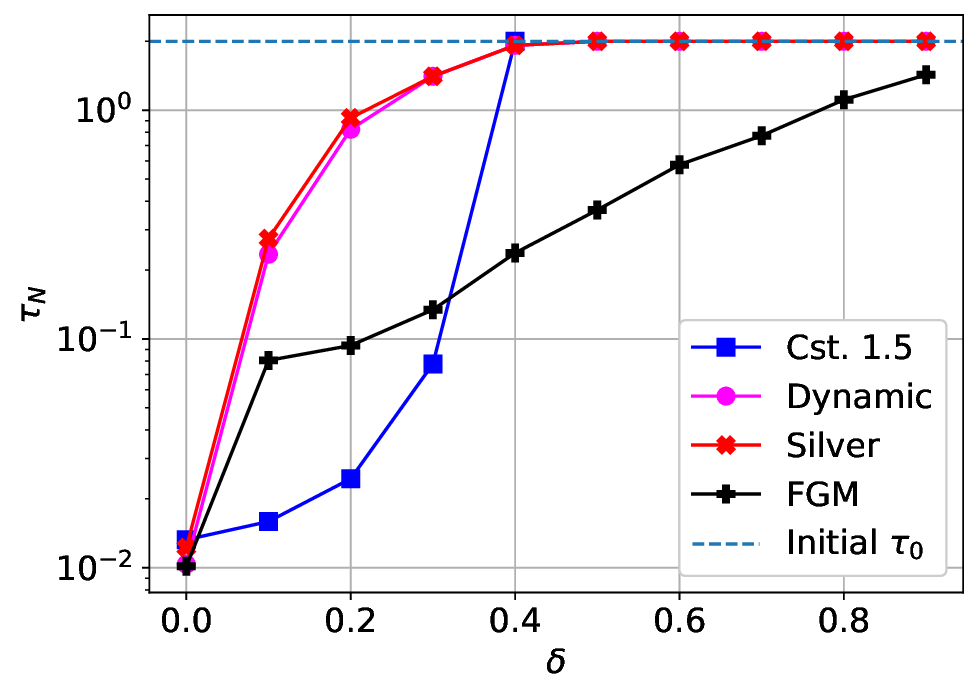}
        \caption{Stepsize schedules from Section \ref{sec::schedules}}
        \label{fig::PEP_exact}
    \end{subfigure}
    \begin{subfigure}{0.49\textwidth}
    \centering
        \includegraphics[scale=0.35]{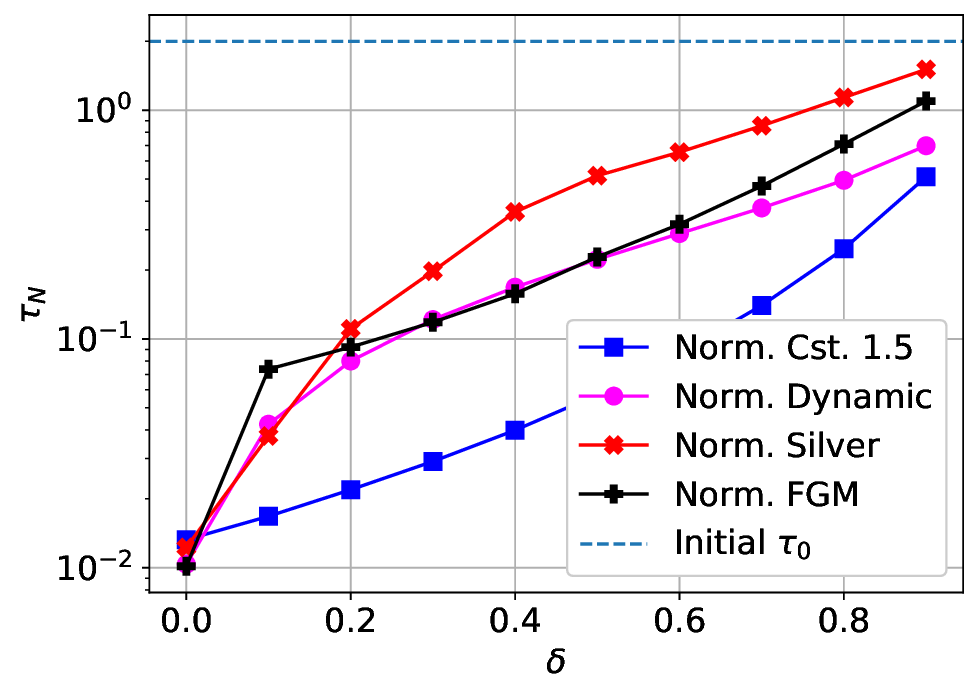}
        \caption{Shortened schedules (divided by $1+\delta$)}
        \label{fig::PEP_normalized}
    \end{subfigure}\vspace*{-.3cm}
    \caption{\centering Exact convergence rate versus inexactness level $\delta$ for $N=50$ iterations}
    \label{fig:PEP}
\end{figure}

A method is said to have a \textit{better} convergence rate than another one if its value $\tau_N$ is \textit{smaller}. Indeed, as $f(x_0)-f(x_*)$ is fixed, the convergence rate indicates the minimum size of the best squared gradient norm among all iterates. In Figure \ref{fig::PEP_exact}, convergence rates of the methods presented in Section \ref{sec::schedules} are presented after $N=50$ iterations for different levels of inexactness. The horizontal line $\tau_0$ represents the value if no iteration is performed (i.e., when staying at $x_0$). \vspace*{-1.0cm}

\vspace*{-0.8cm}

Convergence rates of all methods strongly depend on the inexactness level $\delta$, and no method is particularly robust to the inexactness of the gradient. For example, for an inexactness level $\delta\geq0.4$, most methods do worse than performing no iterations. Indeed, too long step sizes hurt the inexact methods. 
In the context of constant step size gradient descent, it is shown in \cite{vernimmentight} that steps with $h>\tfrac{2}{1+\delta}$ diverge (see Algorithm \ref{algo::inexact_gradient} on the simple univariate function $f(x)=\tfrac{Lx^2}{2}$). We also observed in \cite{vernimmen2024convergence} that the structure of the convergence proof of gradient descent on smooth convex functions can be preserved to tackle inexactness as long as the constant step size satisfies $h\leq\frac{2}{1+\delta}$. Combining both observations naturally leads to the idea that dividing all step sizes by the factor $1+\delta$ can be helpful to tackle inexactness; such schedules will be referred to as \textit{shortened} schedules. Computing the corresponding worst-case rates with PEP 
leads to Figure~\ref{fig::PEP_normalized}, with much more encouraging rates, which suggests using those methods in a practical context. In particular, the shortened constant step size $h=\frac{1.5}{1+\delta}$ appears to be very robust against inexactness, especially for small values of inexactness level $\delta$. The other methods on Figure~\ref{fig::PEP_normalized} are less robust to inexactness, but still benefit from the shortening of the step sizes. 
\vspace{-0.2cm}
\section{Empirical Results}\label{sec::experiments}
\vspace{-0.2cm}
After a theoretical comparison of worst-case convergence rates in the previous section, we now test\footnote{All codes from this section and the previous one can be found on \url{https://github.com/Verpierre/Inexact_PEP_Experiment}} all methods empirically in a more practical setting, both in their original and shortened versions. 

The numerical experiment we present is a \emph{logistic classification problem}, which is formulated as a smooth and convex minimization task, aligned with the theoretical framework introduced in the previous section.

We apply this logistic regression model to a binary image classification task, using a dataset of handwritten images of digits 0 and 1 \cite{arda_mavi_2017}. (We also ran similar tests on another dataset involving phishing websites \cite{website_phishing_379}, and observed similar conclusions.)
Each image is represented as a vector of pixel values $x_k$, and the associated label $y_k \in \{0,1\}$ indicates whether the image shows a zero ($y_k=0$) or a one ($y_k=1$). We denote by $K$ the total number of training images.

The goal is to learn the best classifier—defined by a weight vector $w$ and a bias term $b$—to distinguish between the two digits. This is done by minimizing the logistic loss over all the training examples, leading to the following optimization problem:
\[
\min_{w,b} -\sum_{k=1}^K \left[ y_k \log\left(\frac{1}{1+e^{-(b+w^T x_k)}}\right) + (1-y_k) \log\left(1 - \frac{1}{1+e^{-(b+w^T x_k)}}\right) \right]
\]

This loss function encourages the model to assign high probabilities to the correct labels. Thanks to the properties of the logistic function, this objective is smooth and convex.

For those experiments, the four algorithms presented in Section \ref{sec::schedules}, as well as their shortened versions, are tested with $N=100$ iterations on the logistic classification problem, averaging the results obtained on the same 6 different initial points $x_0$. 
Several levels of inexactness are analyzed, and two main types of relative inexactness are considered: the one from the compressed gradient descent in Section \ref{sec::exp_compressed}, and an adversarial corruption of the gradient (trying to achieve the worst-case behavior) in Section \ref{sec::exp_worst}.

In each plot, and for each value of $\delta$, the best value of the squared norm among all iterates is shown (lower is better), as well as the best accuracy percentage on the training set (higher is better), which may be achieved at different iterates. We also display the corresponding accuracy on the test set (i.e., using the iterate with the best training accuracy). 

The smoothness constant \(L\) of the objective function, required for normalized step sizes, is estimated numerically by tracking the local curvature observed during an optimization process. The procedure begins by selecting an initial parameter \(x_0\) and setting an initial estimate for the smoothness constant, \(L_0 = 0\). We then generate a sequence of iterates \(x_1, \dots, x_N\) using an optimization method that does not require an explicit knowledge of the smoothness constant \(L\). As we progress through the iterates, we incrementally refine our estimate of \(L\) by measuring the local curvature between successive iterates. Specifically, at each step \(k\), the local curvature is calculated as the ratio of the norm of the difference between the gradients at two consecutive iterates to the norm of the difference between the iterates themselves: $L_{\text{local}} = \frac{\|\nabla f(x_k) - \nabla f(x_{k+1})\|}{\|x_k - x_{k+1}\|}$.

We then update our global estimate of the smoothness constant by taking the maximum of the current global estimate and the newly computed local curvature: $L_{k+1} = \max(L_k, L_{\text{local}})$.

This process continues until we have processed all the iterates, at which point the final estimate \(L_N\) represents an upper bound on the local curvature observed throughout the optimization run. In other words, \(L_N\) is the smallest value that satisfies the smoothness inequality for the observed sequence of iterates.

If a subsequent optimization method (such as an inexact gradient method) is now used with the estimated value \(L_N\) and the local curvature observed from these new iterates does not exceed this estimate, then the value of \(L_N\) is likely a good approximation of the true smoothness constant — at least for the set of iterates observed.

\subsection{Empirical analysis of the compressed gradient descent}\label{sec::exp_compressed}
\begin{figure}[H]
    \centering
    \begin{subfigure}{0.49\textwidth}
    \centering
        \includegraphics[scale=0.23]{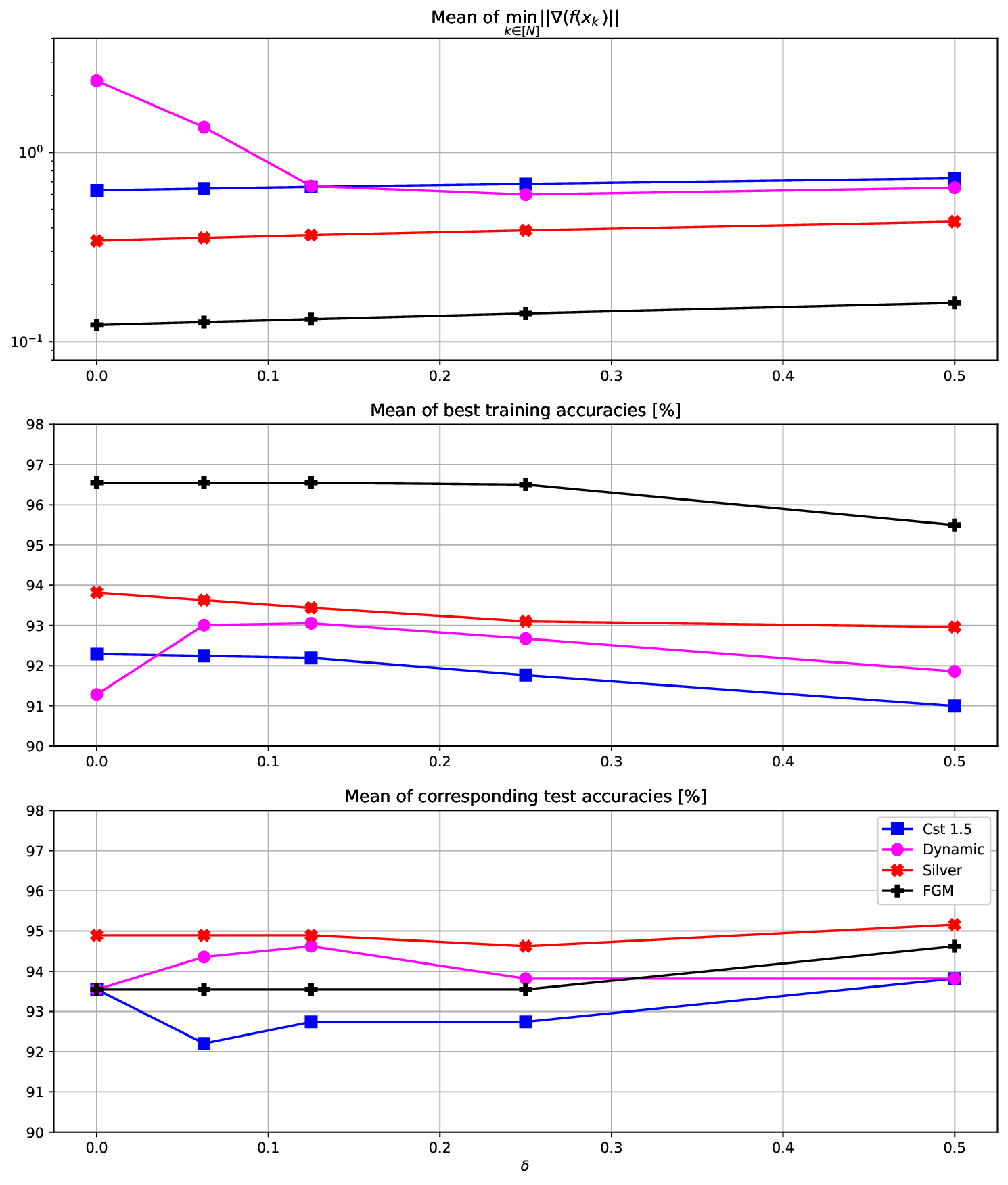}
        \caption{Stepsize schedules from Section \ref{sec::schedules}}
        \label{fig::exp_compressed_exact}
    \end{subfigure}
    \begin{subfigure}{0.49\textwidth}
    \centering
        \includegraphics[scale=0.23]{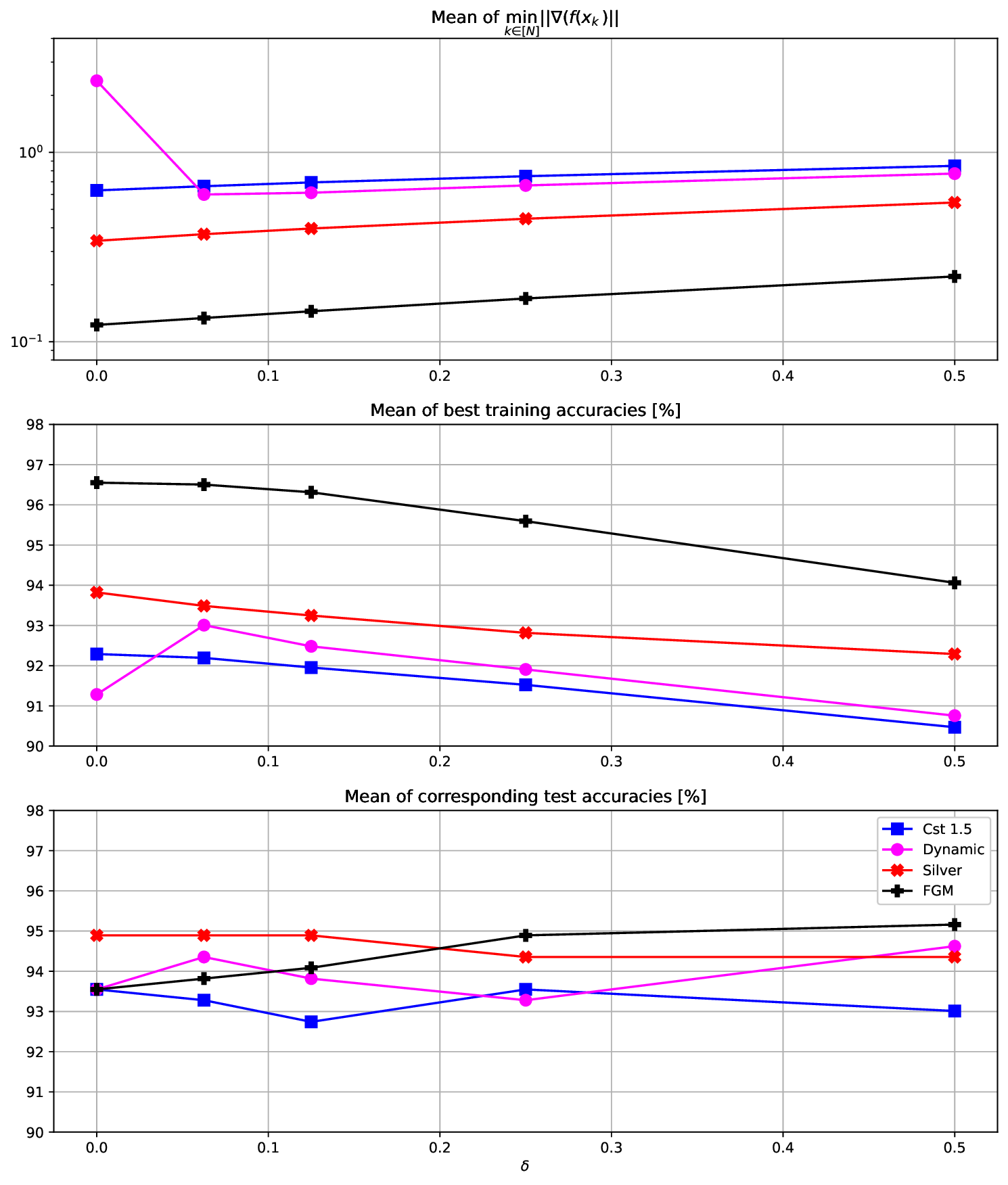}
            \caption{Shortened schedules (divided by $1+\delta$)}
        \label{fig::exp_compressed_normalized}
    \end{subfigure}
    \caption{\centering Empirical results of the compressed gradient descent for different values of inexactness level $\delta_{\text{compressed}}=(\frac{1}{2})^{n_{\text{bit}}+1}$}
    \label{fig::exp_compressed}
\end{figure}

In this first experiment, we consider the compressed gradient descent (Section \ref{sec::truncated_mantissa}) using different numbers of bits $n_{\text{bit}}$ to store the mantissa. In Section \ref{sec::truncated_mantissa}, we saw that the corresponding worst-case value of $\delta$ is equal to $(\frac{1}{2})^{n_{\text{bit}}+1}$, with $n_{\text{bit}}$ being the number of bits used to store each mantissa. The results in this case, for different values of $\delta$, are shown in Figure \ref{fig::exp_compressed}.

\vspace{-0.5cm}

This experiment shows that compressing the mantissa does not seem to cause any significant harm for any method. Shortened schedules also do not seem to make much difference. Furthermore, the ordering of performances in gradient norm and test accuracy between the different methods does not seem to change compared to the exact case. Indeed, FGM beats the Silver Stepsize Schedules method, which in turn beats the other two methods (matching their respective exact complexities $\mathcal{O}(\frac{1}{N^2})$, $\mathcal{O}(\frac{1}{N^{1.27}})$ and $\mathcal{O}(\frac{1}{N})$). It is however interesting to notice that test accuracies for FGM are sometimes worse than for silver step sizes.  In the case of Dynamic Stepsizes, some compression seems to even help the method as $\delta$ grows. This is probably because this compression reduces a bit the size of each component of the gradient, which can equivalently be seen, given the structure of Algorithm \ref{algo::inexact_gradient}, as reducing the step size\footnote{In this sense, it was observed that if the mantissa is truncated by rounding it upwards instead of downwards, the results deteriorated significantly.}. All in all, those results suggest that the error caused by compression is not harmful at all. 
\vspace*{.5cm}

\subsection{Empirical analysis of the gradient descent with corruption error}\label{sec::exp_worst}
\begin{figure}[H]
    \centering
    \begin{subfigure}{0.49\textwidth}
    \centering
        \includegraphics[scale=0.23]{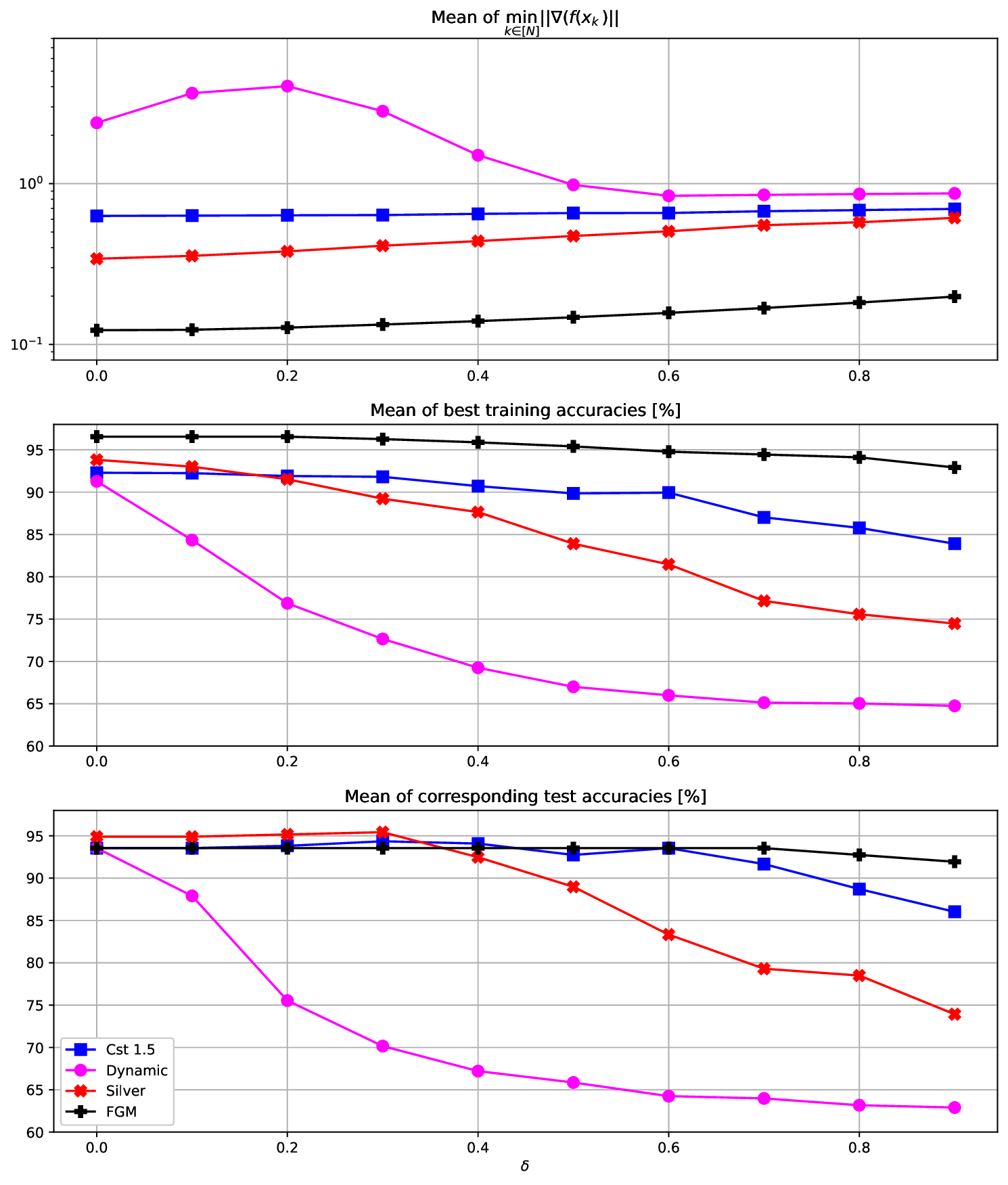}
        \caption{Stepsize schedules from Section \ref{sec::schedules}}
        \label{fig::worst_exact}
    \end{subfigure}
    \begin{subfigure}{0.49\textwidth}
    \centering
        \includegraphics[scale=0.23]{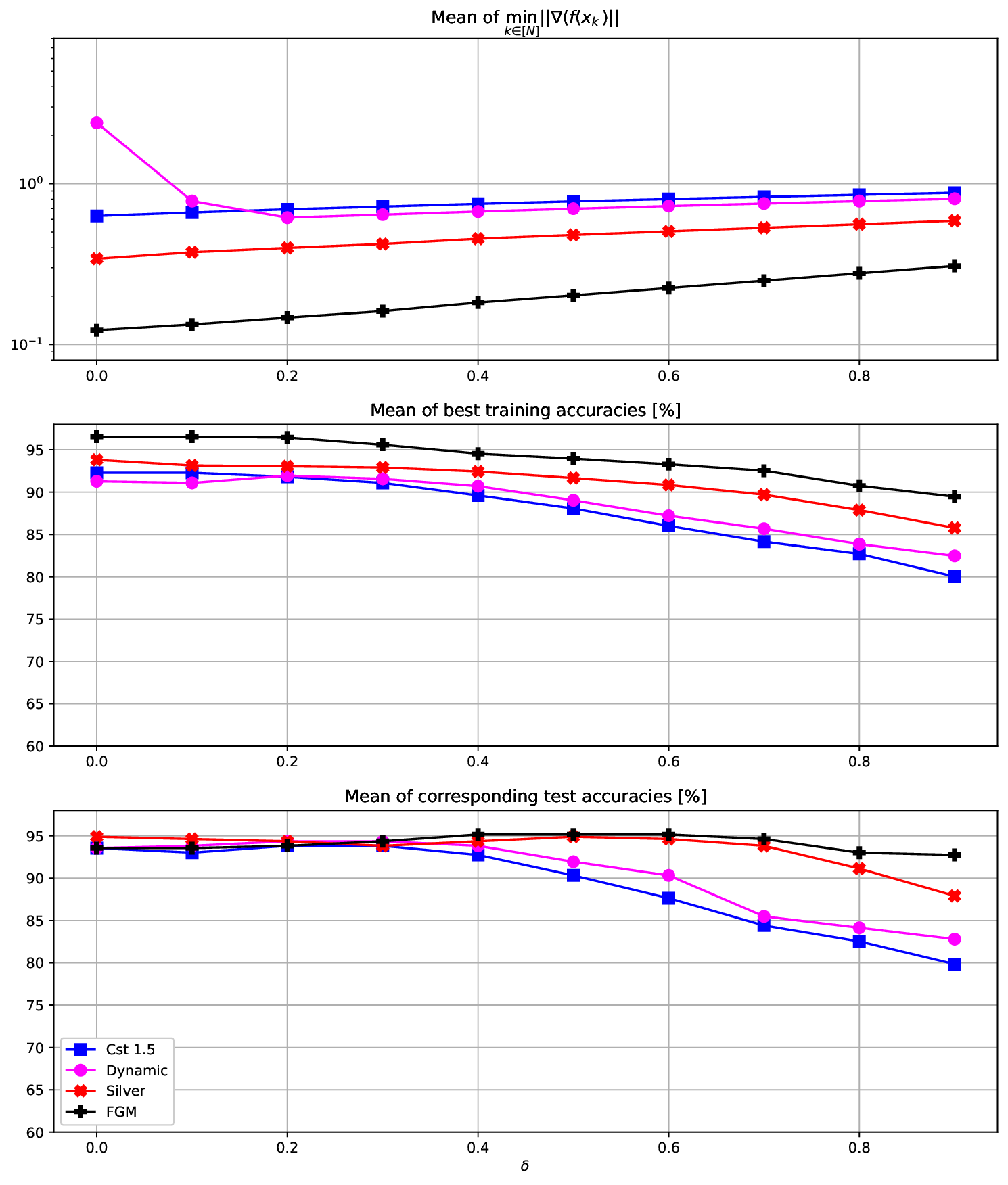}
            \caption{Shortened schedules (divided by $1+\delta$)}
        \label{fig::worst_normalized}
    \end{subfigure}
    \caption{\centering Empirical results of gradient descent with corruption (adversarial) error}
    \label{fig::worst}
\end{figure}

We now consider a more adversarial model of relative inexactness. Specifically, we add to the gradient the longest possible vector (depending on $\delta$) pointing towards the \textit{worst} direction. This \textit{worst} direction is chosen as the direction taking us away from the true minimizer $x_\ast$. In some contexts, for example in cases of corruption error, it is useful to consider a worst-case analysis, as performed in this experiment. Indeed, we can imagine that, in some distributed computing situations, intruders try to make our algorithm fail, see for example \cite{wang2024robust}.
The results given by the experiment for different values of $\delta$ can be seen in Figure \ref{fig::worst} (left without the shortening factor $\frac{1}{1+\delta}$, right with it).
\vspace*{-.5cm}

\vspace{-0.5cm}

In this context, we observe results that are more in line with the theoretical guarantees, as expected. Looking at non-shortened methods, Silver steps and Dynamic steps have more difficulties in handling the inexactness when it is designed to cause them problems, leading to very poor accuracies. Constant steps are also impacted, albeit less significantly. Surprisingly, FGM appears to work well with this kind of adversarial error, even with great inexactness, which is not explained by the worst-case theoretical analysis.

Adding the shortening factor greatly improves the performance of Silver and Dynamic steps, especially for the training and test accuracies (Silver steps being even able to rival FGM for $\delta \le 0.3$). On the other hand, this shortening factor has a small negative effect on the other two methods. Overall, these results show the great performance (and robustness) in practical contexts of both accelerated methods, and we recommend the use of FGM, and of shortened Silver steps when a memoryless method is desired.
\vspace{-0.2cm}
\section{Conclusion}
\vspace{-0.2cm}
We conclude below about the worst-case and practical analysis of the different methods, both without and with the shortening factor $\frac{1}{1+\delta}$. The robustness of the methods is evaluated qualitatively in Table \ref{table::concl}, for \textit{small} inexactness ($\delta<0.4$) and larger ones ($\delta\geq0.4$).
\vspace*{-.8cm}

\begin{table}[ht]
\caption{\centering Qualitative evaluation of robustness based on worst-case and experiments for different $\delta$.
\textcolor{purple}{\protect\ding{55}}=bad, \textcolor{orange}{\protect$\sim$}=intermediate, \textcolor{green}{\protect\checkmark}=good.}
\label{table::concl}
\centering
\begin{tabular}{|c|c|c|c|c|}
\hline
\textbf{Method} & \textbf{$\delta<0.4$} & \textbf{$\delta \geq 0.4$} & \textbf{$\delta<0.4$} & \textbf{$\delta \geq 0.4$} \\
               & \multicolumn{2}{c|}{\textbf{Worst-case}} & \multicolumn{2}{c|}{\textbf{Experiment}} \\
\hline
Cst. 1.5 & \xmark & \xmark & \cmark & \cmark \\
Short. Cst. 1.5 & \cmark & \cmark & \cmark & \cmark \\
\hline
Dynamic & \approxmark & \xmark & \approxmark & \xmark \\
Short. Dynamic & \approxmark & \approxmark & \cmark & \cmark \\
\hline
Silver & \approxmark & \xmark & \cmark & \approxmark \\
Short. Silver & \approxmark & \approxmark & \cmark & \cmark \\
\hline
FGM & \approxmark & \approxmark & \cmark & \cmark \\
Short. FGM & \approxmark & \approxmark & \cmark & \cmark \\
\hline
\end{tabular}
\end{table}

\vspace*{-.5cm}
This shows that, even in the deliberately worst practical scenarios, the methods behave better than what was expected in worst case, which may raise the question of the usefulness of worst-case analysis. However, the worst-case analysis also suggested the use of the shortening factor, which helped all non-accelerated methods in the second set of numerical experiments, still highlighting the relevance of such theoretical worst-case results.
\newline\newline
\noindent \textbf{Future work} Encouraging results from this work open up some possibility of future work. The use of compression of the mantissa could be explored on more complex problems, such as small neural networks with sigmoid activation functions, generalizing the concept of logistic classification. Furthermore, a stochastic analysis of compressed gradient descent, instead of only considering the worst-case compression error, would also be interesting. Indeed, on average, errors caused by compression are much smaller than the worst-possible ones.
\newline\newline
\noindent \textbf{Acknowledgments} The first author is supported by the UCLouvain \emph{Fonds spéciaux de Recherche} (FSR). He would like to thank Nicolas Mil-Homens Cavaco for the fruitful discussions and Hugo Siberdt for the initial code in the exact case on the second dataset. \vspace*{-.5cm}

\bibliographystyle{plain}
\bibliography{my_biblio}

\end{document}